\documentclass[11pt]{amsart}
\usepackage{latexsym,enumerate}
\usepackage{amsmath,amsthm,amsopn,amstext,amscd,amsfonts,amssymb}
\usepackage{color}

\numberwithin{equation}{section}

\numberwithin{theorem}{section}
\numberwithin{corollary}{section}
\numberwithin{lemma}{section}
\numberwithin{definition}{section}
\numberwithin{proposition}{section}
\numberwithin{remark}{section}

\newcommand{\medint}{-\kern  -,375cm\int}

\begin{document}
\title[ A weighted isoperimetric inequality in a wedge ]{ A weighted
isoperimetric inequality in a wedge }
\author{ F. Brock$^1$ - F. Chiacchio$^2$ - A. Mercaldo$^2$}
\thanks{}
\date{}

\begin{abstract}
\noindent Let $c, k_1 , \ldots , k_N $ be non-negative numbers, and define a
measure $\mu $ in the wedge $W:= \{x\in \mathbb{R} ^N :\, x_i >0 , i=1,
\ldots ,N\} $ by $d\mu = e^{c|x|^2 } x_1 ^{k_1 }\cdots x_N ^{k_N } \, dx $.
It is shown that among all measurable subsets of $W$ with fixed $\mu$
-measure, the intersection of $W$ with a ball centered at the origin renders
the weighted perimeter relative to $W$ a minimum. 
\bigskip


\end{abstract}

\maketitle

\setcounter{footnote}{1} \footnotetext{
Leipzig University, Department of Mathematics, Augustusplatz, 04109 Leipzig,
Germany, e-mail: brock@math.uni-leipzig.de}

\setcounter{footnote}{2} \footnotetext{
Dipartimento di Matematica e Applicazioni \textquotedblleft R.
Caccioppoli\textquotedblright , Universit\`{a} degli Studi di Napoli
``Fe\-derico II", Complesso Monte S. Angelo, via Cintia, 80126 Napoli,
Italy, e-mails: francesco.chiacchio@unina.it, mercaldo@unina.it}


\section{Introduction}

Let a measure $\nu $ be defined by $d\nu =\phi (x)\, dx$, where $\phi $ is a
positive Borel measurable function defined on a subset $\Omega $ of
$\mathbb{R} ^N $.
If $M$ is Lebesgue measurable set with $M\subset \Omega $, we define
the $\nu $-measure of $M$
\begin{equation}  \label{intronu}
\nu (M)=\int_{M}d\nu =\int_{M}\phi(x)\, dx
\end{equation}
and the $\nu $-perimeter of $M$ relative to $\Omega$
\begin{equation*}
P_{\nu }(M,\Omega )=\sup \left\{ \int_{M}\mbox{div }(\mathbf{v}(x)\phi(x))\,
dx:\, \mathbf{v}\in C_0 ^1(\Omega ,\mathbb{R}^{N}),\,|\mathbf{v}|\leq 1\
\right\}\, .
\end{equation*}
Note that if $M$ is a smooth set, then
\begin{equation*}
P_{\nu }(M,\Omega )=\int_{\partial M\cap \Omega }\phi(x)\, d{\mathcal{H}}
_{N-1}(x).
\end{equation*}
The isoperimetric problem reads as
\begin{equation}  \label{Imu}
I_{\nu , \Omega }(m) :=\inf \{P_{\nu }(M,\Omega ):\, M\subset \Omega ,\,\nu
(M)=m\},\quad m>0.
\end{equation}
One says that $M$ is an isoperimetric set if $\nu (M)=m$ and $I_{\nu ,
\Omega }(m)= P_{\nu }(M,\Omega )$.
In this paper we consider the case that $\Omega$ is the wedge  $W$ in
$\mathbb{R}^N$, where
\begin{equation}  \label{wedge1}
W:= \{ x \in \mathbb{R}^N :\, x_i >0 , \, i=1 , \ldots ,N \},
\end{equation}
and we  determine functions $\phi$ having $B_{R}\cap W $ as an
isoperimetric set. \newline
Here and throughout the paper, $B_R$ and $B_R(x)$ denote the ball of radius
$R$ centered at zero and at $x$, respectively. \newline
In a recent paper \cite{BCM2}, the authors studied the case that $\Omega $ is the half-space
\begin{equation*}
H:= \{ x\in \mathbb{R}^N :\, x_N >0 \},
\end{equation*}
and, among other things, the following two results
were proved. \\[0.2cm]
\textbf{Theorem A:} \textsl{(see \cite{BCM2}, Theorem 1.1) \newline
Let
\begin{equation*}
\varphi (x) = x_N ^k \mbox{ exp } \{ c |x|^2 \} , \quad x\in H,
\end{equation*}
where $c,k\geq 0 $. Then $I_{\nu ,H } ( m) = P_{\nu } ( B_R \cap H, H) $ for
every $R>0$ such that $m= \nu (B_R \cap H)$.} \\[0.2cm]
\textbf{Theorem B:} \textsl{(see \cite{BCM2}, Theorem 2.1 and Lemma 2.1)
\newline
Let $\phi \in C^2 (W)$, and assume $\phi $ is in separated form,
\begin{equation*}
\phi (x) = \prod _{i=1 } ^N \phi _i (x_i ) ,
\end{equation*}
where $\phi _i \in C^2 ((0, \infty )) \cap C ([0, \infty )) $, $\phi _i (x_i
) >0 $ if $x_i >0 $, ($i=1, \ldots ,N$). Further, suppose that $I_{\nu ,W }
(m) = P_{\nu } (B_R \cap W, W)$ for every $R>0 $ such that $m= \nu (B_R \cap
W) $. Then
\begin{equation}  \label{weight}
\phi (x) = a \mbox{ exp } \{ c|x|^2 \} \prod _{i=1} ^N x_i ^{k_i } , \quad
x\in W,
\end{equation}
where $a>0$, $k_i \geq 0$, ($i=1, \ldots , N$), and $c\in \mathbb{R}$.}
\\[0.2cm]
Theorems A and B are imbedded in a wide bibliography related to the
isoperimetric problems for \textsl{manifolds with density } (see, for
instance,
\cite{BBMP, Bo,  Borell, BCM,  BMP,  CMV, CJQW, DDNT, MS, Mo, Mo2, RCBM, S}).

The following result accomplishes Theorem B, and its proof will be given in
Section 2: \\[0.2cm]
\textbf{Theorem 1 :} 
\textsl{Let $c, k_i $ be nonnegative numbers, ($i= 1, \ldots ,N$), and let a
measure $\mu $ on $W$ be defined by
\begin{equation}  \label{mumeasure}
d\mu := \mbox{ exp } \{ c|x|^2 \} \prod_{i=1 } ^N x_i ^{k_i } \, dx .
\end{equation}
Then $I_{\mu ,W } (m) = P_{\mu } (B_R \cap W, W)$ for every $R>0 $ such that
$m= \mu (B_R \cap W) $. } \\[0.2cm]
 It is customary to write isoperimetric inequalities as a
relation between the perimeter and the measure of a set. \newline
Set
\begin{eqnarray*}
\mathbf{k } & := & (k_1 , \ldots , k_N ), \\
|\mathbf{k}| & := & \sum _{i=1} ^N k_i , \\
h(r) & := & e^{cr ^ 2 } r^{N-1 +|\mathbf{k} | } , \\
H(r) & := & \int_0 ^r e^{ct ^ 2 } t^{N-1 +|\mathbf{k} | } \, dt , \quad
(r\geq 0), \quad \mbox{and } \\
\kappa & := & \int_{{\mathbb{S}} ^{N-1} \cap W } x_1 ^{k_1 } \cdots x_N
^{k_N } \, d{\mathcal{H}} _{N-1} (x).
\end{eqnarray*}
With these notations Theorem 1 can be reformulated as
\begin{equation}  \label{pmu1}
P_{\mu } (M, W) \geq \kappa h\left( H^{-1} ( \mu (M)/\kappa ) \right ) ,
\end{equation}
with equality if $M= B_R \cap W$, for some $R>0 $. Note, in the special case
$c=0$, (\ref{pmu1}) reads
\begin{equation}  \label{pmu2}
P_{\mu } (M, W) \geq \kappa ^{1/(N+|\mathbf{k}|)} \left( ( N+|\mathbf{k}|)
\mu (M) \right) ^{(N-1+|\mathbf{k}|)/(N+|\mathbf{k}|)} .
\end{equation}
Theorem 1 has numerous applications, and they will be analysed in a
forthcoming paper \cite{BCM4}:\newline
The fact that sets $B_R \cap W $, ($R>0$), are isoperimetric for the
weighted measure $\mu $ imply a Polya-Szeg\"{o} - type inequality, comparing
the weighted Sobolev norms of a given function and its weighted
rearrangement (compare with \cite{T}, p. 125). In turn, this allows to find the
best constants in some Sobolev inequalities for functions defined in the
wedge $W$. \newline
Furthermore, Theorem 1 gives rise to sharp comparison result for weighted
elliptic problems in subsets of $W$
(compare with \cite{T0}, \cite{T2}, and \cite{BBMP}).

\section{Proof of the main result}

In this section we prove Theorem 1. A crucial role is played by the
following \\[0.2cm]
\textbf{Lemma 1: } \textsl{Let $k>0$, $l>0$, and define a function $\sigma
\in C^{\infty } ((0, \pi /2 )) \cap C([0, \pi /2 ] )$ implicitly by
\begin{equation}  \label{st}
\int_0 ^{\theta } \sin ^k t \cos ^l t \, dt = c_1 \int _0 ^{\sigma (\theta )
} \sin ^{k+l } s \, ds , \quad \theta \in [0, \pi /2 ],
\end{equation}
where
\begin{equation}  \label{c1}
c_1 := \frac{ \int_0 ^{\pi /2} \sin ^k t \cos ^l t \, dt }{ \int _0 ^{\pi }
\sin ^{k+l } s \, ds } .
\end{equation}
Then }
\begin{equation}  \label{s'>1}
\sigma ^{\prime } (\theta ) \geq 1, \quad \theta \in (0, \pi /2 ).
\end{equation}
\textsl{Proof: \/} Equation (\ref{st}) implies
\begin{equation}  \label{eqn_s'}
\sin ^k \theta \cos ^l \theta = c_1 \sigma ^{\prime } (\theta ) \sin ^{k+l }
\sigma (\theta ) , \quad \theta \in (0, \pi /2) ,
\end{equation}
$\sigma (0)=0$, and $\sigma (\pi /2 )= \pi $. Set
\begin{equation*}
f(\theta ):= \sin ^k \theta \cos ^l \theta -c_1 \sin ^{k+l} \sigma (\theta )
, \quad \theta \in [0, \pi /2 ] ,
\end{equation*}
and note that $f\in C^{\infty } ((0, \pi /2 ))\cap C([0, \pi /2 ]) $, with
$f(0)=f(\pi /2 )=0$. Then (\ref{s'>1}) holds iff
\begin{equation}  \label{f>0}
f(\theta ) \geq 0, \quad \theta \in (0, \pi /2) .
\end{equation}
Assume
(\ref{f>0}) was not true. Then there exists $\theta _0 \in (0, \pi /2)$ with
\begin{equation}  \label{def_t0}
f(\theta _0 )<0, \ f^{\prime }(\theta _0 ) =0, \ f^{\prime \prime }
(\theta_0 )\geq 0.
\end{equation}
This, in turn, also implies that
\begin{equation}  \label{s'<1}
\sigma ^{\prime }(\theta _0 ) <1 .
\end{equation}
By (\ref{eqn_s'}) we have
\begin{equation}  \label{f'simple}
f^{\prime } (\theta ) = c_1 \sigma ^{\prime }(\theta ) \sin ^{k+l } \sigma
(\theta ) \left( k\cot \theta - l\tan \theta - (k+l) \cot \sigma (\theta )
\right) , \quad \theta \in (0, \pi /2).
\end{equation}
The second condition in (\ref{def_t0}) and (\ref{f'simple}) give
\begin{equation}  \label{cond1}
k\cot \theta _0 -l \tan \theta _0 - (k+l ) \cot \sigma (\theta _0 ) =0.
\end{equation}
Then, differentiating (\ref{f'simple}), the third condition in (\ref{def_t0}
), and (\ref{cond1}) give
\begin{equation*}
f^{\prime \prime } (\theta _0 )= c_1 \sigma ^{\prime } (\theta _0 ) \sin
^{k+l } \sigma (\theta _0 ) \left( -\frac{k}{\sin ^2 \theta _0 } - \frac{l}
{\cos ^2 \theta _0 } + \frac{(k+l) \sigma ^{\prime } (\theta _0 ) }{\sin ^2
\sigma (\theta _0 ) } \right) \geq 0,
\end{equation*}
which implies
\begin{equation}  \label{cond3}
1 > \sigma ^{\prime }(\theta _0 ) \geq \left( \frac{k}{\sin ^2 \theta _0 } +
\frac{l}{\cos ^2 \theta _0 } \right) \frac{\sin ^2 \sigma (\theta _0 ) }{k+l}
.
\end{equation}
On the other hand, (\ref{cond1}) yields
\begin{equation*}
\frac{(k \cot \theta _0 - l \tan \theta _0 )^2 }{ (k+l )^2 } = \frac{1-\sin
^2 \sigma (\theta _0 )}{\sin ^2 \sigma (\theta _0 )} ,
\end{equation*}
that is,
\begin{equation}  \label{cond4}
\sin ^2 \sigma (\theta _0 ) = \frac{ (k+l)^2 }{ (k+l)^2 + (k \cot \theta _0
- l \tan \theta _0 )^2 }.
\end{equation}
Plugging (\ref{cond4}) into (\ref{cond3}), we find
\begin{equation}  \label{cond5}
1 > \sigma ^{\prime } (\theta _0 ) \geq \frac{ (k+l ) \left( \frac{k}{\sin
^2 \theta _0 } + \frac{l}{\cos ^2 \theta t_0 } \right) }{ (k+l)^2 + (k\cot
\theta _0 -l \tan \theta _0 ) ^2 }.
\end{equation}
This implies
\begin{equation}  \label{cond6}
(k+l )^2 + (k\cot \theta _0 - l \tan \theta _0 ) ^2 > (k+l) \left( \frac{k}{
\sin ^2 \theta _0 } + \frac{l}{\cos ^2 \theta _0 } \right) .
\end{equation}
Setting $\sin ^2 \theta _0 =: z\in (0,1)$, this yields
\begin{equation*}
(k+l )^2 + k^2 \frac{1-z}{z} + l^2 \frac{z}{1-z} -2kl > (k+l) \left( \frac{k
}{z} +\frac{l}{1-z} \right) ,
\end{equation*}
that is,
\begin{equation*}
\frac{kl}{z } + \frac{kl}{1-z} <0,
\end{equation*}
which is impossible.
Hence (\ref{f>0}) follows, and Lemma 1 is proved.
$\hfill \Box $ \\[0.2cm]

Now we prove Theorem 1 in two steps. We firstly face, using Lemma above, the bidimensional case and then the result is achieved in its full generality via an induction argument over the dimension $N$.

\noindent   \textsl{Proof of Theorem 1. Step 1: The case $N=2$.  }\newline
We write $(x,y)$ for points in $\mathbb{R}^2 $, and $W := \{ (x,y)\in
\mathbb{R}^2 :\, x>0, y>0\}$. Introduce a measure on $W $ by
\begin{equation*}
d\mu = \mbox{exp}\, \{c(x^2 +y^2 ) \} x^l y^k \, dx\, dy,
\end{equation*}
where $c,l$ and $k$ are nonnegative constants. \newline
Let $M$ be a smooth set contained in $W $, and choose $R>0 $ such that $\mu
(B_R \cap W ) = \mu (M)$. Then
\begin{equation}  \label{decomp1}
\partial M \cap W = \bigcup_{k=1} ^m C_k ,
\end{equation}
where the $C_k $'s are mutually non-intersecting, smooth curves, and
\begin{equation}  \label{perimeter1}
P_{\mu } (M, W ) = \sum_{k=1} ^m \int _{C_k } \mbox{exp}\, \{ c(x^2 + y^2
)\} x^l y^k \, ds , \quad (ds: \ \mbox{ Euclidean arc length differential }).
\end{equation}
Let $C$ be one of the curves in the decomposition (\ref{decomp1}), with
parametrization
\begin{equation*}
C:= \{ (x(t),y(t)): \, t\in [a,b]\} ,
\end{equation*}
where $x,y\in C^1 ( [a,b])$ and $a,b\in \mathbb{R} $, $a<b$. Then
\begin{equation}  \label{perimeter2}
\int _{C } \mbox{exp}\, \{ c(x^2 + y^2 )\} x^l y^k \, ds = \int_a ^b
\mbox{exp}\, \{ c(x^2 (t)+y^2 (t) ) \} x^l (t) y^k (t) \sqrt{(x^{\prime }
(t))^2 + (y^{\prime } (t) )^2 } \, dt .
\end{equation}
Let $p: W\to (0, +\infty )\times (0, \pi /2) $ be the polar coordinates
mapping given by
\begin{equation*}
p (r\cos \theta , r\sin \theta ) := (r, \theta ), \quad (r> 0, \, 0<\theta <
\pi /2 ),
\end{equation*}
and write $(r(t), \theta (t)):= p (x(t), y(t)) $, ($a\leq t\leq b$). Then
also
\begin{equation}  \label{perimeter3}
\int _{C } \mbox{exp}\, \{ c(x^2 + y^2 )\} x^l y^k \, ds = \int_a ^b e^{cr^2
(t) } r^{k+l } (t) \sqrt{ r^2 (t) (\theta ^{\prime } (t))^2 + (r^{\prime }
(t) )^2 } \, \sin ^k \theta (t) \cos ^l \theta (t) \, dt .
\end{equation}
(Here and in the following, "$\, ^{\prime } \, $" denotes differentiation
w.r.t. $t$.) \newline
Next, let $U: (0, +\infty )\times (0, \pi /2 ) \to (0, +\infty ) \times (0,
\pi )$ be the mapping given by
\begin{equation*}
U(r, \theta ) := (r, \sigma (\theta )) , \quad (r>0, \, 0<\theta <\pi /2 ),
\end{equation*}
where $\sigma $ is the function defined by (\ref{st}) and (\ref{c1}).
Finally, let $\mathbb{R} ^2 _+ $ be the upper half-plane $\{ (u , v )\in
\mathbb{R} ^2 :\, v >0 \} $, and define a diffeomorphism from $W$ onto
$\mathbb{R} ^2 _+ $ by
\begin{equation*}
T:= p^{-1} \circ U\circ p.
\end{equation*}
Writing $(u, v )$ for points in $\mathbb{R} ^2 _+ $, and
$(x,y)\equiv (r \cos \theta , r \sin \theta )$ for points in $W$, we have
\begin{equation*}
(u,v) \equiv T(x,y) = (r\cos \sigma (\theta ), r \sin \sigma (\theta ) ) .
\end{equation*}
Introduce a measure $\widetilde{\mu } $ on $\mathbb{R} ^2 _+ $ by
\begin{equation*}
d\widetilde{\mu} := \mbox{exp}\, \{ c(u ^2 + v ^2 ) \} v ^{k+l} \, du \, dv
.
\end{equation*}
Using the notation $(u(t),v(t)) := T(x(t),y(t))$, and since $p (u(t), v(t))
= (r(t), \sigma (\theta ((t)))$, ($a\leq t\leq b$), we find, similarly as
above,
\begin{equation}  \label{perimeter4}
P_{\widetilde{\mu} } (T(M), \mathbb{R} ^2 _+ ) = \sum _{k=1} ^m
\int_{T(C_k )} \mbox{exp}\, \{ c(u ^2 + v ^2 ) \} v ^{k+l} \, ds ,
\end{equation}
and
\begin{eqnarray} \label{perimeter5}
& & \int_{T(C )} \mbox{exp}\, \{ c(u ^2 + v ^2 ) \} v ^{k+l} \, ds \\
& = & \int_a ^b \mbox{exp}\, \{ c(u ^2 (t) + v ^2 (t) ) \} v ^{k+l} (t)
\sqrt{ (u^{\prime } (t))^2 + (v^{\prime } (t)) ^2 } \, dt  \notag \\
& = & \int_a ^b e^{cr^2 (t) } r^{k+l} (t) \sqrt{ r^2 (t) \left( \frac{d\sigma
}{d\theta} \right) ^2 (\theta ^{\prime } (t))^2 + (r ^{\prime } (t)) ^2 } \,
\sin ^{k+l } \sigma (\theta (t)) \, dt .  \notag
\end{eqnarray}
By Lemma 1 we have $\frac{d\sigma }{d \theta } \geq 1$, which implies,
together with (\ref{perimeter5}), (\ref{st}) and (\ref{c1}),
\begin{eqnarray}  \label{perimeter6}
& & \int_{T(C )} \mbox{exp}\, \{ c(u ^2 + v ^2 ) \} v ^{k+l} \, ds \\
& \leq & \int_a ^b e^{cr^2 (t) } r^{k+l} (t) \sqrt{ r^2 (t) (\theta ^{\prime
} (t))^2 + (r ^{\prime } (t)) ^2 } \frac{d\sigma }{d\theta } \, \sin ^{k+l}
\sigma (\theta (t) ) \, dt  \notag \\
& = & \frac{1}{c_1 } \int_a ^b e^{cr^2 (t) } r^{k+l} (t) \sqrt{ r^2 (t)
(\theta ^{\prime } (t))^2 + (r ^{\prime } (t)) ^2 } \, \sin ^k \theta (t)
\cos ^l \theta (t) \, dt  \notag \\
& = & \frac{1}{c_1 } \int _C \exp \, \{ c( x^2 + y^2 ) \} x^l y^k \, ds ,
\notag
\end{eqnarray}
with equality if $r(t)=$ const. for $t\in [a,b]$. In view of
(\ref{perimeter1}), (\ref{perimeter3}), (\ref{perimeter4}) and (\ref{perimeter6})
we conclude that
\begin{equation}  \label{perimeter7}
c_1 P_{\widetilde{\mu} }( T(M), \mathbb{R} ^2 _+ ) \leq P_{\mu } (M, W) ,
\end{equation}
with equality if $M= B_R \cap W$. \newline
Finally, using (\ref{st}), an elementary calculation shows that
\begin{eqnarray}  \label{measure1}
\mu (M) & = & \int\limits _M \! \! \int \mbox{exp}\, \{ c(x^2 + y^2 )\} x^l
y^k \, dx\, dy \\
& = & \int\limits _{ p(M)} \! \! \! \! \int e^{cr^2 } r^{k+l+1 } \sin ^k
\theta \cos ^l \theta \, dr\, d\theta  \notag \\
& = & c_1 \int\limits _{ p(M)} \! \! \! \! \int e^{cr^2 } r^{k+l+1 } \sin
^{k+l} \sigma (\theta ) \frac{d\sigma }{d\theta } \, dr\, d\theta  \notag \\
& = & c_1 \int\limits _{p(T(M))} \! \! \! \! \int e ^{cr^2 } r^{k+l+1 } \sin
^{k+l} \sigma \, dr\, d\sigma  \notag \\
& = & c_1 \widetilde{\mu} (T(M)).  \notag
\end{eqnarray}
Now Theorem A, for $N=2$, tells us that
$P_{\widetilde{\mu} } (T(M), \mathbb{R} ^2 _+ ) \geq P_{\widetilde{\mu} } (T(B_R \cap W ), \mathbb{R} ^2 _+ )$.
Together with (\ref{perimeter7}) and (\ref{measure1}) this proves the
assertion for smooth sets. \newline
If $M$ is a measurable subset of $W$ with finite $\mu$-perimeter, then, by
the very properties of the weighted perimeter, there exists a sequence of
smooth sets $\{ M_n \} $, ($M_n \subset W$, $n\in \mathbb{N}$), such that
$\lim_{n\to \infty } \mu (M_n \Delta M ) =0$, and $\lim_{n\to \infty } P_{\mu
} (M_n ) = P_{\mu } (M)$. Hence the assertion for $M$ follows by
approximation. $\hfill \Box $ \\[0.2cm]
\textsl{Step 2: The general case. } \newline
We proceed by induction over the dimension $N$. We write $y= (x ^{\prime
},x_N , x_{N+1 } ) $ for points in $\mathbb{R}^{N+1} $, where $x^{\prime
}\in \mathbb{R}^{N} $, and $x_{N+1 } \in \mathbb{R} $, and
\begin{equation*}
W_{N+1} := \{ y= (x^{\prime }, x_N , x_{N+1} ) \in \mathbb{R} ^{N+1} : \,
x_i >0 , \, i=1,\ldots ,N+1 \} .
\end{equation*}
Assume that the assertion holds true for some $N\in \mathbb{N} $,
($N\geq 2 $).
\newline
Let a measure $\mu $ on $W_{N+1} $ be given by
\begin{equation*}
d \mu = \mbox{ exp } \{ c (|x^{\prime}|^2 + x_N ^2 + x_{N+1} ^2 ) \} \,
\prod_{i=1 } ^{N+1} x_i ^{k_i }\, dy .
\end{equation*}
We define two measures $\nu _1 $ and $\nu _2 $ by
\begin{eqnarray*}
d \nu _1 & = & \mbox{ exp } \{ c |x^{\prime}|^2  \}\, \prod_{i=1 } ^{N-1} x_i
^{k_i } \, dx^{\prime }, \\
d \nu _2 & = & x_{N} ^{k_N } x_{N+1} ^{k_{N+1}} \mbox{ exp } \{ c ( x_N ^2 +
x_{N+1} ^2 ) \} \, dx_N dx_{N+1 } ,
\end{eqnarray*}
and note that $d\mu = d\nu _1 d\nu _2 $. \newline
Let $M $ be a subset of $W_{N+1}$ having finite and positive $\mu $-measure.
\newline
We define $2$-dimensional slices
\begin{equation*}
M(x^{\prime }) := \{ (x_N , x_{N+1 } ):\, (x^{\prime }, x_N , x_{N+1} ) \in
M\} , \quad (x^{\prime }\in \mathbb{R}^{N-1} ).
\end{equation*}
Let $M^{\prime }:= \{ x^{\prime }\in \mathbb{R}^{N-1} :\, 0< \mu _2
(M(x^{\prime })) \} $, and note that $\nu _2 (M(x^{\prime })) < +\infty $
for a.e. $x^{\prime }\in M^{\prime }$. For all those $x^{\prime }$, let
$Q(x^{\prime })$
be the quarter disc in $\mathbb{R} ^2 _+ $ centered at
$(0,0)$ with $\nu_2 (M(x^{\prime })) = \nu _2 (Q(x^{\prime }))$. (For
convenience, we put $Q(x^{\prime })= \emptyset $ for all $x^{\prime }\in
M^{\prime }$ with $\nu _2 (M(x^{\prime })) =+\infty $.) Let
\begin{equation*}
W_2 := \{ (x_N , x_{N+1} ):\, x_N >0, \, x_{N+1 } >0 \} .
\end{equation*}
Since Theorem 1 holds in the two-dimensional case, we have that
\begin{equation}  \label{isopnu2}
P_{\nu_2 } (Q(x^{\prime }), W_2 ) \leq P_{\nu _2 } (M(x^{\prime }), W_2 )
\quad \mbox{for a.e. } \ x^{\prime }\in M^{\prime }.
\end{equation}
Let
\begin{equation*}
Q:= \{ y=(x^{\prime }, x_N , x_{N+1} ) :\, (x_N , x_{N+1} )\in Q(x^{\prime
}) ,\, x^{\prime }\in M^{\prime }\} .
\end{equation*}
The product structure of the measure $\mu $ tells us that \newline
\textbf{(i)} $\mu (M)= \mu (Q)$, and \newline
\textbf{(ii)} the isoperimetric property for slices, (\ref{isopnu2}),
carries over to $M$, that is,
\begin{equation}  \label{isopnu}
P_{\mu } (Q, W_{N+1 }) \leq P_{\mu } (M, W_{N+1} ) ,
\end{equation}
(see for instance Theorem 4.2 of \cite{BBMP}). \newline
Note again, the slice $Q (x^{\prime }) = \{ ( x_N , x_{N+1 } ) :\,
(x^{\prime }, x_N , x_{N+1 } ) \in Q \} $ is a quarter disc $\{ (r\cos
\theta , r \sin \theta ):\, 0 < r < R(x^{\prime }) ,\, \theta \in (0, \pi/2
) \} $, with $0<R(x^{\prime }) <+\infty $, ( $x^{\prime }\in M^{\prime }$).
Set
\begin{equation*}
K:= \{ (x^{\prime }, r) : \, 0< r <  R(x^{\prime }), \, x^{\prime }\in
M^{\prime }\} ,
\end{equation*}
and introduce a measure $\alpha $ on
\begin{equation*}
W_N := \{ (x^{\prime },r ) :\, x_i >0 , \, i=1,\ldots ,N-1, \, r>0 \}
\end{equation*}
by
\begin{equation*}
d\alpha := a r^{k_N+k_{N+1}+1 } \mbox{ exp } \{ c(|x^{\prime}|^2+ r^2 )\} \,
dx^{\prime }dr,
\end{equation*}
where
\begin{equation*}
a := \int_0 ^{\pi / 2} \mbox{cos} ^{k_N} \theta \,\,\mbox{sin} ^{k_{N+1}
}\theta \, d\theta.
\end{equation*}
An elementary calculation then shows that
\begin{equation*}
\mu (Q) = \alpha (K),
\end{equation*}
and
\begin{equation*}
P_{\mu } (Q, W_{N+1} ) = P_{\alpha } (K, W_N ).
\end{equation*}
Let $B_R $ denote the open ball in $\mathbb{R}^N $ centered at the origin,
with radius $R$, and choose $R>0 $ such that
\begin{equation*}
\alpha (B_R \cap W_N ) = \alpha ( K).
\end{equation*}
By the induction assumption it follows that
\begin{equation}  \label{isopalpha}
P_{\alpha } (B_R \cap W_N , W_N ) \leq P_{\alpha } (K, W_N ).
\end{equation}
Finally, let $M^{\bigstar } $ defined by
\begin{equation*}
M^{\bigstar } := \{ y= (x^{\prime }, x_N , x_{N+1 } ):\, |x^{\prime}|^2+ x_N
^2 + x_{N+1 } ^2 < R^2 , \, x_i >0 ,\, i=1,..,N+1 \} .
\end{equation*}
Then

\begin{equation*}
\mu (M^{\bigstar }) = \mu (M)
\end{equation*}
and
\begin{equation*}
P_{ \mu } (M^{\bigstar} , W_{N+1} ) = P_{\alpha } (B_R \cap W_N , W_N ).
\end{equation*}

The equalities above, together with (\ref{isopalpha}) and (\ref{isopnu}), yield
\begin{equation*}
\label{isopN+1} P_{\mu } (M^{\bigstar }, W_{N+1} ) \leq P_{\mu } (M, W_{N+1}
),
\end{equation*}
that is, the isoperimetric property holds for $N+1 $ in place of $N$
dimensions. The Theorem is proved. $\hfill \Box $


\begin{thebibliography}{99}






\bibitem{BBMP} \textsc{M.F. Betta, F. Brock, A. Mercaldo, M.R. Posteraro},
Weighted isoperimetric inequalities on $\mathbb{R}^{n}$ and applications to
rearrangements. \textsl{Math. Nachr.} \textbf{281} (2008), no. 4, 466--498.



\bibitem{Bo} \textsc{C. Borell}, The Brunn-Minkowski inequality in Gauss
space. \textsl{Invent. Math.} \textbf{30} (1975), no. 2, 207--211.

\bibitem{Borell} \textsc{C. Borell}, The Ornstein-Uhlenbeck velocity process
in backward time and isoperimetry. (1986), Preprint Chalmers University of
Technology 1986-03/ISSN 0347-2809.

\bibitem{BCM} {\sc F. Brock, F. Chiacchio, A. Mercaldo}, A class of
degenerate elliptic equations and a Dido's problem
with respect to a measure.
{\sl J. Math. Anal. Appl.} {\bf 348} (2008), no. 1, 356--365.

\bibitem{BCM2} \textsc{F. Brock, F. Chiacchio, A. Mercaldo}, Weighted
isoperimetric inequalities in cones and applications. \textsl{Nonlinear
Analysis T.M.A.} \textbf{75} (2012), no. 15, 5737–-5755. 

\bibitem{BCM4} \textsc{F. Brock, F. Chiacchio, A. Mercaldo}, A weighted
isoperimetric inequality in a wedge and applications, \textsl{in preparation}.

\bibitem{BMP} \textsc{F. Brock, A. Mercaldo, M.R. Posteraro}, On Schwarz and
Steiner symmetrization with respect to a measure. to appear in : \textsl{
Revista Matem\'{a}tica Iberoamericana}.




\bibitem{CMV} \textsc{A. Ca\~{n}ete, M. Miranda Jr., D. Vittone}, Some
isoperimetric problems in planes with density. \textsl{J. Geom. Anal.}
\textbf{20} (2010), no.2, 243--290.

\bibitem{CJQW} \textsc{C. Carroll, A. Jacob, C. Quinn, R. Walters}, The
isoperimetric problem on planes with density. \textsl{Bull. Aust. Math. Soc.}
\textbf{78} (2008), no. 2, 177--197.


\bibitem{DDNT} \textsc{J. Dahlberg, A. Dubbs, E. Newkirk, H. Tran},
Isoperimetric regions in the plane with density $r^p $. \textsl{New York J.
Math.} \textbf{16} (2010), 31--51.









\bibitem{MS} \textsc{C. Maderna, S. Salsa}, Sharp estimates of solutions to
a certain type of singular elliptic boundary value problems in two
dimensions. \textsl{Applicable Anal.} \textbf{12} (1981), no. 4, 307--321.





\bibitem{Mo} \textsc{F. Morgan}, Manifolds with density. \textsl{Notices
Amer. Math. Soc.} \textbf{52} (2005), no.8, 853--858.

\bibitem{Mo2} \textsc{F. Morgan}, The Log-Convex Density Conjecture. 
\textsl{Contemporary Mathematics} \textbf{545} (2011), 209--211.





\bibitem{RCBM} \textsc{C. Rosales, A. Ca\~{n}ete, V. Bayle, F. Morgan}, On
the isoperimetric problem in Euclidean space with density. \textsl{Calc.
Var. Partial Differential Equations} \textbf{31} (2008), no. 1, 27--46.


\bibitem{S} \textsc{V.N. Sudakov, B.S. Cirel'son}, Extremal properties of
half-spaces for spherically invariant measures. \textsl{Zap. Nau\v cn. Sem.
Leningrad. Otdel. Mat. Inst. Steklov. (LOMI)} \textbf{41} (1974), 14--24.

\bibitem{T0} \textsc{G. Talenti}, Soluzioni a simmetria assiale di equazioni
ellittiche. \textsl{Ann. Mat. Pura Appl.} \textbf{73}(4) (1966) 127--158.

\bibitem{T2} \textsc{G. Talenti}, Elliptic equations and rearrangements.
\textsl{Ann. Sc. Norm. Sup. Pisa Cl. Sci.} \textbf{3} (1976), no. 4,
697--718.

\bibitem{T} \textsc{G. Talenti}, A weighted version of a rearrangement
inequality. \textsl{Ann. Univ. Ferrara, Sez. VII (N.S.)} \textbf{43} (1997),
121--133.


\end{thebibliography}
\end{document}